\documentclass[12pt]{article}
\usepackage{tikz}
\usetikzlibrary{calc,through,backgrounds}
\usepackage{multicol}

\usepackage{bm,amssymb,amsthm,amsfonts,amsmath,latexsym,cite,color, psfrag,graphicx,ifpdf,pgfkeys,pgfopts,xcolor,extarrows}

\newtheorem{theo}{Theorem}[section]

\newtheorem{lem} [theo]{Lemma}
\newtheorem{coro}[theo]{Corollary}

\allowdisplaybreaks

\makeatletter \@addtoreset{equation}{section}
\@addtoreset{theo}{}\makeatother

\usepackage{hyperref}

\def\qed{\hfill \rule{4pt}{7pt}}
\def\pf{\noindent {\it Proof.} }
\def\S{  \mathfrak{S}}

\def\GL{   \mathrm{GL}  }
\def\Char{   \mathrm{char}  }
\def\tr{   \mathrm{tr}  }
\def\F{   \mathcal{F}  }
\def\Newton{  \mathrm{Newton}}
\def\wt{  \mathrm{wt}}
\def\opt{  \mathrm{opt}}
\begin{document}
\begin{center}
{\Large\bf Lattice Points in the Newton Polytopes of

Key Polynomials}

\vskip 4mm
{\footnotesize  NEIL  J.Y. FAN, PETER L. GUO, SIMON C.Y. PENG, SOPHIE C.C. SUN }

%
%
%
%
\end{center}

\noindent{A{\scriptsize BSTRACT}.}
We confirm a conjecture of Monical, Tokcan and Yong on
a characterization of the lattice points in
the Newton polytopes of key polynomials.

\section{Introduction}

Key polynomials  $\kappa_\alpha(x)$ associated to compositions $\alpha\in \mathbb{Z}_{\geq 0}^n$,
also called   Demazure characters,
 are characters of the Demazure
modules for the general linear groups \cite{Dem-1,Dem-2}. They are nonsymmetric polynomial
generalizations of Schur polynomials. Key polynomials are intimately connected  with other important
polynomials in algebraic combinatorics. For example, every Schubert polynomial is a positive sum
 of key polynomials (see for example Lascoux and  Sch\"utzenberger \cite{Las-2}, Reiner and Shimozono \cite{Rei}),
 every key polynomial is a positive sum of
 Demazure atoms (see for example Haglund,  Luoto,  Mason  and  van Willigenburg\cite{Hag},
  Lascoux and  Sch\"utzenberger \cite{Las-1},   Mason \cite{Mas}). Moreover,  $\kappa_\alpha(x)$ can be
 realized as a   specialization of
the nonsymmetric Macdonald polynomial $E_\alpha(x; q, t)$,
that is, $\kappa_\alpha(x)=E_\alpha(x; q=\infty, t=\infty)$, see Ion \cite{Ion}.

This paper is concerned with the   Newton
polytope of $\kappa_\alpha(x)$.
Given a polynomial
\[f=\sum_{\alpha\in \mathbb{Z}_{\geq 0}^n}c_\alpha x^\alpha \in \mathbb{R}[x_1,\ldots,x_n],\]
  the Newton polytope  of $f$ is the convex hull of the exponent vectors of $f$:
\[\Newton(f)=\mathrm{conv}(\{\alpha\colon c_\alpha\neq 0\}).\]
By definition, each exponent vector of $f$ is a lattice point in $\Newton(f)$.
If every lattice point in $\Newton(f)$
is also an  exponent vector of $f$, then we say that $f$ has saturated Newton polytope (SNP).
The SNPness  of polynomials has been investigated by Monical, Tokcan and Yong \cite{Mon}.
They \cite[Conjecture 3.10]{Mon} conjectured that key polynomials have the SNPness property.
This conjecture was confirmed by
Fink,   M\'esz\'aros  and   St.$\,$Dizier \cite{Fin}.
It was also conjectured by Monical, Tokcan and Yong \cite[Conjecture 3.13]{Mon} and proven
 by Fan and Guo \cite{Fan} that the vertices of $\Newton(\kappa_\alpha)$ can be generated
 by permutations in a lower  interval in the Bruhat order.

 Monical, Tokcan and Yong \cite[Conjecture 3.11]{Mon} further conjectured a characterization
   of the
 lattice points in $\Newton(\kappa_\alpha)$, or equivalently, of the exponent vectors
 of $\kappa_\alpha(x)$. The task of this paper is to prove this conjecture.
Let $\alpha=(\alpha_1,\ldots, \alpha_n)\in \mathbb{Z}_{\geq 0}^n$.
For $1\leq i<j\leq n$, let $t_{i,j}(\alpha)$ be the composition obtained from $\alpha$
by interchanging $\alpha_i$ and $\alpha_j$, and let
\[m_{i,j}(\alpha)=\alpha+e_i-e_j,\]
 where $e_k\, (1\leq k\leq n)$ is the standard coordinate vector.
 For a vector  $\beta\in \mathbb{Z}_{\geq 0}^n$, define $\beta\leq_\kappa \alpha$
 if $\beta$ can be generated from $\alpha$ by applying a sequence of moves $t_{i,j}$ for $\alpha_i<\alpha_j$, and  $m_{i,j}$ for $\alpha_i<\alpha_j-1$.

 \begin{theo}[\mdseries{Monical-Tokcan-Yong \cite[Conjecture 3.11]{Mon}}]\label{conj}
A vector $\beta$ is a lattice point in $\Newton(\kappa_\alpha)$ if and only
if $\beta\leq_\kappa \alpha$. Equivalently, $\beta$ is an exponent vector of
$\kappa_\alpha(x)$ if and only if $\beta\leq_\kappa \alpha$.
 \end{theo}

When the parts of $\alpha$ are weakly increasing,   $\kappa_\alpha(x)$
equals the Schur polynomial $s_{\lambda}(x)$, where $\lambda$ is the partition obtained by
rearranging the parts of $\alpha$ decreasingly, see for example Reiner and Shimozono \cite{Rei}.
In this case, the Newton polytope of $s_{\lambda}(x)$ is
$\mathcal{P}_{\lambda}$, the permutohedron whose vertices are   rearrangements of $\lambda$.
A classical theorem due to Rado \cite{Rad} states that for a partition $\mu$,
$\mathcal{P}_{\mu}\subseteq \mathcal{P}_{\lambda}$ if and only if
 $\mu\unlhd \lambda$ in the dominance order.
By Rado's theorem, it is easy to check
that the lattice points of $\mathcal{P}_{\lambda}$ are the
rearrangements of partitions $\mu$  with $\mu\unlhd \lambda$.
On the other hand, it can be shown that   when  $\alpha$ is weakly increasing,
$\beta\leq_\kappa \alpha$ if and only if $\beta$ is a  rearrangement
of some partition $\mu\unlhd \lambda$.
This yields a proof of Theorem
\ref{conj} for  the case when $\alpha$ is weakly increasing.

The following corollary is a direct consequence of Theorem \ref{conj},
which generalizes Rado's theorem from partitions to compositions.

\begin{coro}
The Newton polytope $\Newton(\kappa_\beta)$ of $\kappa_\beta(x)$
is contained in the Newton polytope $\Newton(\kappa_\alpha)$ of $\kappa_\alpha(x)$
if and only if $\beta\leq_\kappa \alpha$.
\end{coro}

We also remark that Theorem \ref{conj} leads to a description of the lattice
points in certain Bruhat interval polytopes.
For
two permutations $u\leq v$ in the Bruhat order, the    Bruhat interval
polytope $\textsf{Q}_{u,v}$  is the convex hull of the permutations in
the Bruhat interval $[u, v]$.  Bruhat interval polytopes were introduced
by Kodama and Williams \cite{Kod}, and their combinatorial properties were  studied
by Tsukerman and Williams \cite{Tsu}. When $\alpha=w$ is a permutation of $[n]=\{1,2,\ldots,n\}$,
that is, the parts of $\alpha$ are rearrangements of $1,2,\ldots,n$,
the Newton polytope $\Newton(\kappa_w)$ coincides with $\textsf{Q}_{w,w_0}$
\cite[Corollary 1.3]{Fan},
where $w_0=n\cdots 2 1$. Hence,  when $\alpha=w$ is a permutation of $[n]$,
 Theorem \ref{conj} characterizes  the lattice
points in the Bruhat interval polytope $\textsf{Q}_{w,w_0}$.

We prove Theorem \ref{conj} for any compositions by employing
   the realization of key polynomials as the
dual characters of   flagged Weyl modules associated to skyline diagrams,
see Section \ref{Weyl}.
The structure of flagged Weyl modules has been used by
Fink,   M\'esz\'aros  and   St.$\,$Dizier \cite{Fin} to
prove the SNPness  of Schubert polynomials and
key polynomials. In Section 3, using  the flagged Weyl module associated to the
skyline diagram $D(\alpha)$ of $\alpha$, we encode the monomials
appearing in $\kappa_\alpha(x)$ in terms of column-strict flagged
fillings of $D(\alpha)$.
We introduce   operations on column-strict flagged
fillings  of $D(\alpha)$,  which enable us to reflect the moves $t_{i,j}$
and $m_{i,j}$ from compositions to column-strict flagged
fillings of $D(\alpha)$. This allows us to obtain
 a proof of Theorem \ref{conj}.

\section{Key polynomials and   flagged Weyl modules}\label{Weyl}



Key polynomials  can be  defined using the Demazure operators $\pi_i=\partial_ix_i$.
Here, $\partial_i$ is the  divided difference operator, that is,
given a polynomial $f(x)\in \mathbb{Z}[x_1,\ldots,x_n]$,
 $\partial_i$ sends $f(x)$ to
\[\partial_i (f(x))=\frac{f(x)
    -s_i f(x)}{x_i-x_{i+1}},\]
where $s_i f(x)$ is obtained from $f(x)$ by exchanging $x_i$ and $x_{i+1}$.

If $\alpha$ is a partition (that is, the parts of $\alpha$ are weakly
decreasing), then set
$\kappa_\alpha(x)=x^\alpha.$
Otherwise, choose $i$ such that $\alpha_i<\alpha_{i+1}$. Let $\alpha'$ be the
composition obtained from $\alpha$ by interchanging $\alpha_i$ and $\alpha_{i+1}$,
namely, $\alpha'=t_{i,i+1}(\alpha)$. Set
\[\kappa_\alpha(x)=\pi_i (\kappa_{\alpha'}(x))  =\partial_i(x_i\kappa_{\alpha'}(x)).\]
The above definition is independent of the choice of
the position $i$,  since the Demazure operators satisfy the
  braid relations: $\pi_i \pi_j=\pi_j \pi_i$ for $|i-j|>1$, and
$\pi_i \pi_{i+1} \pi_i= \pi_{i+1}\pi_i  \pi_{i+1}$.
For example, for $\alpha=(1,3,2)$, we have
\begin{align*}
\kappa_{(1,3,2)}(x)&=\pi_1\, \kappa_{(3,1,2)}(x)=\pi_1 \pi_2\,\kappa_{(3,2,1)}(x)=\pi_1 \pi_2\, (x_1^3x_2^2x_3)\\[5pt]
&=x_1^3x_2^2x_3+x_1^3x_2x_3^2+x_1^2x_2^3x_3+x_1^2x_2^2x_3^2+x_1x_2^3x_3^2.
\end{align*}

In the remaining of this section, we    briefly review the structure of
flagged Weyl modules associated to diagrams of an $n\times n$ grid.
In particular, the key polynomial $\kappa_\alpha(x)$ is equal to the dual character of
the flagged Weyl module  associated to the skyline diagram of $\alpha$.
The flagged Weyl modules
 can be  constructed  by means  of determinants \cite{Mag}.
Here we follow the notation in \cite{Fin}.

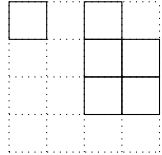
\begin{figure}[h]
\begin{center}
\begin{tikzpicture}

\def\rectanglepath{-- +(5mm,0mm) -- +(5mm,5mm) -- +(0mm,5mm) -- cycle}

\draw [step=5mm,dotted] (0mm,0mm) grid (20mm,20mm);
\draw (0mm,15mm) \rectanglepath;\draw (10mm,15mm) \rectanglepath;
\draw (10mm,10mm) \rectanglepath;\draw (15mm,10mm) \rectanglepath;
\draw (10mm,5mm) \rectanglepath;\draw (15mm,5mm) \rectanglepath;

\end{tikzpicture}
\end{center}
\vspace{-6mm}
\caption{A diagram of a $4\times 4$ grid.}
\label{diagram}
\end{figure}

A diagram $D$ is  a collection of boxes of
an $n\times n$ grid. We use $(i,j)$ to denote the box of the $n\times n$ grid  in row $i$
and column $j$, where  the row indices increase from top to bottom and the column indices increase  from left to right.
With this notation, a diagram $D$
can be written as an ordered list $D=(D_1,D_2,\ldots, D_n)$ of $n$ subsets of $[n]$,
 this is, $i\in D_j$ if and only of    $(i,j)$ is a box of $D$.  For example,
the diagram in  Figure \ref{diagram} can be represented as $(\{1\}, \emptyset, \{1,2,3\},\{2,3\})$.

Let $\GL(n,\mathbb{C})$ be the group of $n\times n$ invertible matrices over $\mathbb{C}$, and
let $B$ be the subgroup   consisting of the $n\times n$ upper-triangular matrices.
Let $Y$ be the $n\times n$ upper-triangular matrix whose entries are
 indeterminates $y_{ij}$ where $i\leq j$. Denote by $\mathbb{C}[Y]$
 the ring of polynomials in the variables $\{y_{ij}\}_{ i\leq j}$.
The group $\GL(n,\mathbb{C})$ acts on $\mathbb{C}[Y]$ (on the right)
as follows.  Given a matrix $g\in \GL(n,\mathbb{C})$ and a polynomial $f(Y)\in \mathbb{C}[Y]$, define
\[f(Y)\cdot g=f(g^{-1} Y).\]
For two diagrams $C=(C_1,\ldots,C_n)$ and $D=(D_1,\ldots,D_n)$, write  $C\leq D$ if
$C_j\leq D_j$ for every $1\leq j\leq n$, where $C_j\leq D_j$ means that $|C_j|=|D_j|$ and
for $1\leq k\leq |C_j|$,
the $k$-th least element of $C_j$ is less than or equal to the $k$-th least element of $D_j$.
The flagged Weyl module $\mathcal{M}_D$ associated to a diagram $D$ is
a $B$-module  defined by
\begin{equation}\label{MM}
\mathcal{M}_D= \mathrm{Span}_\mathbb{C}\left\{\prod_{j=1}^n \det\left(Y_{D_j}^{C_j}\right)\colon C\leq D \right\},
\end{equation}
where, for two subsets $R$ and $S$ of $[n]$, $Y^R_S$ denotes the submatrix
of $Y$ with row indices in $R$ and column indices in $S$.
It should be noted  that $\prod_{j=1}^n \det\left(Y_{D_j}^{C_j}\right)\neq 0$ if and only if $C\leq D$.

Let $X=\mathrm{diag}(x_1,\ldots,x_n)$ be a  diagonal matrix, which can
be viewed as a linear transformation from $\mathcal{M}_D$ to $\mathcal{M}_D$ via the
$B$-action. The character of $\mathcal{M}_D$ is defined  as the trace of $X$:
\[\Char(\mathcal{M}_D )(x)=\tr(X\colon \mathcal{M}_D\rightarrow \mathcal{M}_D).\]
The dual character of $\mathcal{M}_D$ is the character of the dual module $\mathcal{M}_D^*$, which is given by
\begin{align*}
\Char^\ast(\mathcal{M}_D)(x)&=
\tr(X\colon \mathcal{M}_D^*\rightarrow \mathcal{M}_D^*)\\[5pt]
&=\Char(\mathcal{M}_D )(x_1^{-1},\ldots, x_n^{-1}).
\end{align*}

Two families of flagged Weyl modules are of particular interest.
The first one is the flagged Weyl module associated to the Rothe diagram
$D(w)$ of a permutation
$w$ of $[n]$. In this case,
Kra\'skiewicz and   Pragacz \cite{Kra-1,Kra-2} showed that
the Schubert polynomial $\S_w(x)$ of $w$ is equal to the dual character
of   $\mathcal{M}_{D(w)}$. The second one is the   flagged Weyl module
associated to the skyline diagram of a composition $\alpha$, which is
the structure that we need for the purpose of this paper.

\begin{figure}[h]
\begin{center}
\begin{tikzpicture}

\def\rectanglepath{-- +(5mm,0mm) -- +(5mm,5mm) -- +(0mm,5mm) -- cycle}

\draw [step=5mm,dotted] (50mm,0mm) grid (70mm,20mm);
\draw (50mm,15mm) \rectanglepath;
\draw (50mm,10mm) \rectanglepath;\draw (55mm,10mm) \rectanglepath; \draw (50mm,0mm) \rectanglepath;
\draw[dotted](50mm,0mm)--(50mm,20mm);

\end{tikzpicture}
\end{center}
\vspace{-6mm}
\caption{ The skyline diagram of  $\alpha=(1,2,0,1)$.}
\label{RS}
\end{figure}
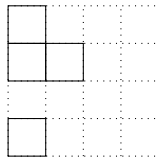

The skyline
diagram  $D(\alpha)$ of a composition $\alpha$
 is the diagram consisting of the first $\alpha_i$ boxes in row $i$.
For example,
Figure \ref{RS} illustrates the skyline diagram of $\alpha=(1,2,0,1)$.

\begin{theo}[\mdseries{Demazure \cite{Dem-2}}]\label{BB}
Let $D(\alpha)$ be the skyline
diagram of a composition $\alpha$. Then
\begin{equation}
\kappa_\alpha(x)=\Char^\ast(\mathcal{M}_{D(\alpha)})(x).
\end{equation}
\end{theo}

Notice that, for $C\leq D$, the effect of the action of the diagonal matrix
$X=\mathrm{diag}(x_1,\ldots,x_n)$ on the polynomial $\prod_{j=1}^n \det\left(Y_{D_j}^{C_j}\right)$ is
\[\prod_{j=1}^n \det\left(Y_{D_j}^{C_j}\right)\cdot X= \prod_{j=1}^n\prod_{i\in C_j}x_i^{-1}\cdot \prod_{j=1}^n \det\left(Y_{D_j}^{C_j}\right).\]
This implies that  the polynomial $\prod_{j=1}^n \det\left(Y_{D_j}^{C_j}\right)$  is an
eigenvector of $X$ with eigenvalue
\[\prod_{j=1}^n\prod_{i\in C_j}x_i^{-1}.\]
For any diagram $C=(C_1,\ldots,C_n)$, let
\[x^C=\prod_{j=1}^n\prod_{i\in C_j}x_i.\]
Therefore, the set of monomials appearing in the dual character $\Char^\ast(\mathcal{M}_{D})(x)$
is exactly
$\left\{x^C\colon C\le D\right\}.$
Restricting $D$ to a skyline diagram $D(\alpha)$ and combining Theorem \ref{BB},
we are given the following description of   monomials
appearing in $\kappa_\alpha(x)$.

\begin{theo}\label{KK}
The set of monomials appearing in $\kappa_\alpha(x)$ is
\[\left\{x^C\colon C\le D(\alpha)\right\}.\]
\end{theo}

\section{Proof of Theorem \ref{conj}}\label{Sec3}

To provide a proof of Theorem \ref{conj}, we first characterize the   monomials appearing in
$\kappa_\alpha(x)$   in terms of  certain fillings
of  skyline diagrams.

\begin{figure}[h]
\begin{center}
\begin{tabular}{ccc}
\begin{tikzpicture}[scale=1.2]
\def\rectanglepath{-- +(4mm,0mm) -- +(4mm,4mm) -- +(0mm,4mm) -- cycle}
\draw [step=4mm,dotted] (0mm,0mm)grid (20mm,20mm);
\draw (0mm,8mm) \rectanglepath;
\draw (0mm,4mm) \rectanglepath;
\draw (0mm,0mm) \rectanglepath;
\draw (4mm,12mm) \rectanglepath;
\draw (4mm,8mm) \rectanglepath;
\draw (4mm,4mm) \rectanglepath;

\draw (8mm,12mm) \rectanglepath;
\draw (8mm,4mm) \rectanglepath;
\draw (8mm,0mm) \rectanglepath;
\draw (12mm,4mm) \rectanglepath;
\draw (12mm,0mm) \rectanglepath;
\draw (16mm,16mm) \rectanglepath;
\draw (16mm,8mm) \rectanglepath;
\draw (16mm,4mm) \rectanglepath;

\node at (2mm,10mm) {\small{2}};
\node at (2mm,6mm) {\small{1}};
\node at (2mm,2mm) {\small{4}};

\node at (6mm,14mm) {\small{1}};
\node at (6mm,10mm) {\small{3}};
\node at (6mm,6mm) {\small{2}};

\node at (10mm,14mm) {\small{1}};
\node at (10mm,6mm) {\small{2}};
\node at (10mm,2mm) {\small{4}};

\node at (14mm,6mm) {\small{4}};
\node at (14mm,2mm) {\small{2}};

\node at (18mm,18mm) {\small{1}};
\node at (18mm,10mm) {\small{3}};
\node at (18mm,6mm) {\small{2}};

\node at (10mm,-6mm) {(a)};
\end{tikzpicture}&
\hspace{1.5cm}
\begin{tikzpicture}[scale=1.2]
\def\rectanglepath{-- +(4mm,0mm) -- +(4mm,4mm) -- +(0mm,4mm) -- cycle}
\draw [step=4mm,dotted] (0mm,0mm)grid (20mm,20mm);
\draw (0mm,8mm) \rectanglepath;
\draw (0mm,4mm) \rectanglepath;
\draw (0mm,0mm) \rectanglepath;
\draw (4mm,12mm) \rectanglepath;
\draw (4mm,8mm) \rectanglepath;
\draw (4mm,4mm) \rectanglepath;

\draw (8mm,12mm) \rectanglepath;
\draw (8mm,4mm) \rectanglepath;
\draw (8mm,0mm) \rectanglepath;
\draw (12mm,4mm) \rectanglepath;
\draw (12mm,0mm) \rectanglepath;
\draw (16mm,16mm) \rectanglepath;
\draw (16mm,8mm) \rectanglepath;
\draw (16mm,4mm) \rectanglepath;

\node at (2mm,10mm) {\small{1}};
\node at (2mm,6mm) {\small{2}};
\node at (2mm,2mm) {\small{4}};

\node at (6mm,14mm) {\small{1}};
\node at (6mm,10mm) {\small{2}};
\node at (6mm,6mm) {\small{3}};

\node at (10mm,14mm) {\small{1}};
\node at (10mm,6mm) {\small{2}};
\node at (10mm,2mm) {\small{4}};

\node at (14mm,6mm) {\small{2}};
\node at (14mm,2mm) {\small{4}};

\node at (18mm,18mm) {\small{1}};
\node at (18mm,10mm) {\small{2}};
\node at (18mm,6mm) {\small{3}};

\node at (10mm,-6mm) {(b)};
\end{tikzpicture}
\end{tabular}
\caption{Two fillings of a diagram $D$.}\label{FFFsf}
\end{center}
\end{figure}
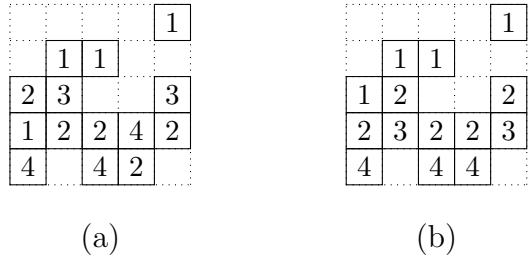

Let $D=(D_1,D_2,\ldots, D_n)$ be a diagram of $[n]^2$.
A filling $F$ of  $D$ is an assignment of positive integers into
the boxes of $D$.
A filling $F$ is called column-strict if the integers  in each column of
$F$ are distinct, and $F$ is called  flagged if for each box in row $i$,
the integer assigned in it does not
 exceed $i$. We denote by $\F(D)$  the set of column-strict flagged
fillings  of $D$.
We also define $\F_{\leq}(D)$ to be the subset of $\F(D)$ consisting of the fillings
$F\in \F(D)$ such that the integers in each column of $F$ are increasing from top to bottom.
For example, Figure \ref{FFFsf}(a) is a filling in $\F(D)$, while Figure \ref{FFFsf}(b)
is a filling in $\F_{\leq}(D)$.
For a filling $F\in \F(D)$, write $\wt(F)=(v_1,v_2,\ldots, v_n)$ to be the weight of $F$,  where $v_i$ is the number of appearances of
$i$ in $F$.

\begin{theo}\label{CCC}
The set of monomials appearing in $\kappa_\alpha(x)$ is
\[\left\{x^{\wt(F)}\colon F\in \F(D(\alpha))\right\}.\]
\end{theo}

\pf For any diagram $D$, there is a direct bijection between $\F_{\leq}(D)$ and the set $\{C\colon C\leq D\}$.
To be specific, given a filling $F\in \F_{\leq}(D)$, let $C=(C_1, C_2,\ldots, C_n)$ be the diagram
such that for $1\leq j\leq n$, $C_j$ is the set of integers filled in the $j$-th
column of $F$. Therefore, by Theorem \ref{KK}, the set of monomials appearing in $\kappa_\alpha(x)$ is
$\left\{x^{\wt(F)}\colon F\in \F_{\leq}(D(\alpha))\right\}$. We conclude the proof by showing that for any diagram $D$,
\begin{equation}\label{AA}
\left\{x^{\wt(F)}\colon F\in \F_{\leq}(D)\right\}=
\left\{x^{\wt(F)}\colon F\in \F(D)\right\}.
\end{equation}
It follows from $\F_{\leq}(D)\subseteq\F(D)$ that
$
\left\{x^{\wt(F)}\colon F\in \F_{\leq}(D)\right\} \subseteq
\left\{x^{\wt(F)}\colon F\in \F(D)\right\}.
$

Now we verify the reverse inclusion. Given a filling $F\in\F(D)$, let $F'$ be a filling obtained from $F$ by resorting the integers in each column increasingly from top to bottom. For example, if
$F$ is the filling in Figure \ref{FFFsf}(a), then   $F'$ is the
filling in Figure \ref{FFFsf}(b).
 By \cite[Proposition 3.5]{Fan},
$F'$ belongs to $\F_{\leq}(D)$. Since $F$ and $F'$ have the same weight,
we see that
\[
\left\{x^{\wt(F)}\colon F\in \F(D)\right\} \subseteq
\left\{x^{\wt(F)}\colon F\in \F_{\leq}(D)\right\}.
\]
This proves \eqref{AA}, and so the proof is complete.
\qed


We now prove the necessity of Theorem \ref{conj}.

\begin{theo}\label{NN}
If $\beta$ is an exponent vector of
$\kappa_\alpha(x)$, then $\beta\leq_\kappa \alpha$.
\end{theo}

To prove Theorem \ref{NN}, we  introduce an operation, called   {\it optimization},  on the fillings
of $\F(D)$. Given a filling $F\in\F(D)$, the optimization of $F$, denoted   $\opt(F)$, is obtained  by rearranging the integers in each column of $F$ as follows.
 Write $D=(D_1,\ldots, D_n)$.
For $1\leq m\leq n$, let $C_m$
denote the set of integers filled in the  $m$-th column $F_m$ of $F$.
Suppose that
\[C_m\cap D_m=\{i_1<i_2<\cdots<i_k\}.\]
Let us construct the $m$-th column $\overline{F}_m$ of $\opt(F)$.
Set $\overline{F}^{(0)}_m=F_m$.
For $1\leq r\leq k$, $\overline{F}_m^{(r)}$ is obtained from $\overline{F}_m^{(r-1)}$ as below. If $i_r$ is filled in the box $(i_r,m)$ of $\overline{F}_m^{(r-1)}$, then let $\overline{F}_m^{(r)}=\overline{F}_m^{(r-1)}$.
Otherwise,  $i_r$ is filled in a box $(t,m)$ of $\overline{F}_m^{(r-1)}$ with $t>i_r$.
Let $\overline{F}_m^{(r)}$
be obtained from $\overline{F}_m^{(r-1)}$  by interchanging $i_r$ and the integer
filled in the box  $(t,m)$.
Define $\overline{F}_m=\overline{F}_m^{(k)}$. By the above construction, each
column $\overline{F}_m^{(r)}$ satisfies the flag constraint. Hence
 $\opt(F)$  is a filling belonging to $\F(D)$. Moreover,  $\opt(F)$  has  the same weight as $F$.

\begin{figure}[h]
\begin{center}
\begin{tabular}{ccc}
\begin{tikzpicture}[scale=1.2]
\def\rectanglepath{-- +(4mm,0mm) -- +(4mm,4mm) -- +(0mm,4mm) -- cycle}
\draw [step=4mm,dotted] (0mm,0mm)grid (20mm,20mm);
\draw (0mm,8mm) \rectanglepath;
\draw (0mm,4mm) \rectanglepath;
\draw (0mm,0mm) \rectanglepath;
\draw (4mm,12mm) \rectanglepath;
\draw (4mm,8mm) \rectanglepath;
\draw (4mm,4mm) \rectanglepath;

\draw (8mm,12mm) \rectanglepath;
\draw (8mm,4mm) \rectanglepath;
\draw (8mm,0mm) \rectanglepath;
\draw (12mm,4mm) \rectanglepath;
\draw (12mm,0mm) \rectanglepath;
\draw (16mm,16mm) \rectanglepath;
\draw (16mm,8mm) \rectanglepath;
\draw (16mm,4mm) \rectanglepath;

\node at (2mm,10mm) {\small{2}};
\node at (2mm,6mm) {\small{1}};
\node at (2mm,2mm) {\small{4}};

\node at (6mm,14mm) {\small{1}};
\node at (6mm,10mm) {\small{3}};
\node at (6mm,6mm) {\small{2}};

\node at (10mm,14mm) {\small{1}};
\node at (10mm,6mm) {\small{2}};
\node at (10mm,2mm) {\small{4}};

\node at (14mm,6mm) {\small{4}};
\node at (14mm,2mm) {\small{2}};

\node at (18mm,18mm) {\small{1}};
\node at (18mm,10mm) {\small{3}};
\node at (18mm,6mm) {\small{2}};

\node at (36mm,10mm){$\xlongrightarrow[]{\text{Optimization}}$};
\end{tikzpicture}&

\begin{tikzpicture}[scale=1.2]
\def\rectanglepath{-- +(4mm,0mm) -- +(4mm,4mm) -- +(0mm,4mm) -- cycle}
\draw [step=4mm,dotted] (0mm,0mm)grid (20mm,20mm);
\draw (0mm,8mm) \rectanglepath;
\draw (0mm,4mm) \rectanglepath;
\draw (0mm,0mm) \rectanglepath;
\draw (4mm,12mm) \rectanglepath;
\draw (4mm,8mm) \rectanglepath;
\draw (4mm,4mm) \rectanglepath;

\draw (8mm,12mm) \rectanglepath;
\draw (8mm,4mm) \rectanglepath;
\draw (8mm,0mm) \rectanglepath;
\draw (12mm,4mm) \rectanglepath;
\draw (12mm,0mm) \rectanglepath;
\draw (16mm,16mm) \rectanglepath;
\draw (16mm,8mm) \rectanglepath;
\draw (16mm,4mm) \rectanglepath;

\node at (2mm,10mm) {\small{2}};
\node at (2mm,6mm) {\small{4}};
\node at (2mm,2mm) {\small{1}};

\node at (6mm,14mm) {\small{2}};
\node at (6mm,10mm) {\small{3}};
\node at (6mm,6mm) {\small{1}};

\node at (10mm,14mm) {\small{2}};
\node at (10mm,6mm) {\small{4}};
\node at (10mm,2mm) {\small{1}};

\node at (14mm,6mm) {\small{4}};
\node at (14mm,2mm) {\small{2}};

\node at (18mm,18mm) {\small{1}};
\node at (18mm,10mm) {\small{3}};
\node at (18mm,6mm) {\small{2}};
\end{tikzpicture}
\end{tabular}
\caption{$\opt(F)$ for $F$ being the filling in Figure \ref{FFFsf}(a).}\label{OO}
\end{center}
\end{figure}
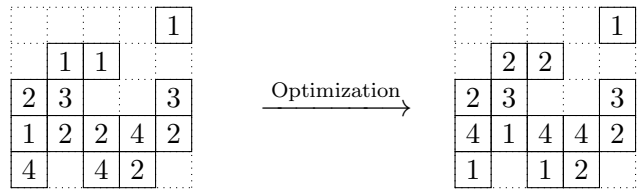

For example, if $F$
is the filling in Figure \ref{FFFsf}(a), then $\opt(F)$ is
the filling as given in Figure \ref{OO}.

%
%
%

Using the optimization operation, we have the following lemma,
which is crucial to the proof of Theorem \ref{NN}.

\begin{lem}\label{LL}
Let $\alpha=(\alpha_1,\ldots,\alpha_n)$ be a composition, and
  $F$ be a filling in  $\F(D(\alpha))$ with $\wt(F)\neq \alpha$. Then there exists a
filling $F'\in\F(D(\alpha))$   such that $\wt(F)=t_{i,j}(\wt(F'))$  or
$\wt(F)=m_{i,j}(\wt(F'))$.
\end{lem}

\pf  We aim to construct a filling $F'\in\F(D(\alpha))$ from
the optimization $\opt(F)$ of $F$ such that
$\wt(\opt(F))=t_{i,j}(\wt(F'))$ or $\wt(\opt(F))=m_{i,j}(\wt(F')).$

Since $\wt(\opt(F))=\wt(F)\neq \alpha$, there must exist a box of $\opt(F)$ that is filled with
an integer not equal to its row index.
Locate the topmost row of $\opt(F)$,
say row $j$,
which contains an entry not equal to its row index.
Suppose that $i$ is the leftmost entry in the $j$-th row of $\opt(F)$
that is not equal to $j$. Since $\opt(F)$  satisfies the flag condition,
 we have $i<j$. Moreover, let $(j,h)$ be the box in row $j$ that contains this
  leftmost entry $i$. Write $\wt(\opt(F))=(\beta_1,\ldots, \beta_n)$.
 The construction of $F'$
 depends on the relative order of $\beta_i$ and $\beta_j$.

Case 1. $\beta_i\le \beta_j$. The filling $F'$ is obtained from
$\opt(F)$ by replacing the entry $i$ in the box $(j, h)$ with $j$.
We explain that $F'$ is a filling in $\F(D(\alpha))$.
Clearly, $F'$ satisfies the flag condition. We still need to verify that
$F'$ is column strict. By the construction of $\opt(F)$,   there is no integer
in the $h$-th column of $\opt(F)$ that is equal to $j$, since otherwise the box $(j,h)$
of $\opt(F)$ would be filled with $j$. This implies that $F'$ is column strict,
and thus $F'\in \F(D(\alpha))$. By the construction of $F'$, we see that
\[\beta=m_{i,j}(\wt(F')).\]

\begin{figure}[htbp]
\begin{center}
\begin{tabular}{ccc}

\begin{tikzpicture}
\def\rectanglepath{-- +(4mm,0mm) -- +(4mm,4mm) -- +(0mm,4mm) -- cycle}
\draw [step=4mm,dotted] (0mm,0mm)grid (32mm,32mm);
\draw (0mm,24mm) \rectanglepath;
\draw (0mm,20mm) \rectanglepath;
\draw (0mm,12mm) \rectanglepath;
\draw (0mm,8mm) \rectanglepath;
\draw (0mm,4mm) \rectanglepath;
\draw (0mm,0mm) \rectanglepath;
\draw (4mm,24mm) \rectanglepath;
\draw (4mm,12mm) \rectanglepath;
\draw (4mm,8mm) \rectanglepath;
\draw (4mm,0mm) \rectanglepath;
\draw (8mm,24mm) \rectanglepath;
\draw (8mm,8mm) \rectanglepath;
\draw (8mm,0mm) \rectanglepath;
\draw (12mm,24mm) \rectanglepath;
\draw (12mm,8mm) \rectanglepath;
\draw (16mm,8mm) \rectanglepath;
\draw (20mm,8mm) \rectanglepath;

\node at (2mm,26mm) {\small{1}};
\node at (2mm,22mm) {\small{3}};
\node at (2mm,14mm) {\small{2}};
\node at (2mm,10mm) {\small{5}};
\node at (2mm,6mm) {\small{4}};
\node at (2mm,2mm) {\small{6}};

\node at (6mm,26mm) {\small{2}};
\node at (6mm,14mm) {\small{5}};
\node at (6mm,10mm) {\small{6}};
\node at (6mm,2mm) {\small{7}};

\node at (10mm,26mm) {\small{2}};
\node at (10mm,10mm) {\small{5}};
\node at (10mm,2mm) {\small{6}};

\node at (14mm,26mm) {\small{2}};
\node at (14mm,10mm) {\small{6}};

\node at (18mm,10mm) {\small{3}};

\node at (22mm,10mm) {\small{3}};
\node at (46mm,14mm){$\xlongrightarrow[]{\text{Optimization}}$};
\end{tikzpicture}&

\begin{tikzpicture}
\def\rectanglepath{-- +(4mm,0mm) -- +(4mm,4mm) -- +(0mm,4mm) -- cycle}
\draw [step=4mm,dotted] (0mm,0mm)grid (32mm,32mm);
\draw (0mm,24mm) \rectanglepath;
\draw (0mm,20mm) \rectanglepath;
\draw (0mm,12mm) \rectanglepath;
\draw (0mm,8mm) \rectanglepath;
\draw (0mm,4mm) \rectanglepath;
\draw (0mm,0mm) \rectanglepath;
\draw (4mm,24mm) \rectanglepath;
\draw (4mm,12mm) \rectanglepath;
\draw (4mm,8mm) \rectanglepath;
\draw (4mm,0mm) \rectanglepath;
\draw (8mm,24mm) \rectanglepath;
\draw (8mm,8mm) \rectanglepath;
\draw (8mm,0mm) \rectanglepath;
\draw (12mm,24mm) \rectanglepath;
\draw (12mm,8mm) \rectanglepath;
\draw (16mm,8mm) \rectanglepath;
\draw (20mm,8mm) \rectanglepath;

\node at (2mm,26mm) {\small{2}};
\node at (2mm,22mm) {\small{3}};
\node at (2mm,14mm) {\small{5}};
\node at (2mm,10mm) {\small{6}};
\node at (2mm,6mm) {\small{4}};
\node at (2mm,2mm) {\small{1}};

\node at (6mm,26mm) {\small{2}};
\node at (6mm,14mm) {\small{5}};
\node at (6mm,10mm) {\small{6}};
\node at (6mm,2mm) {\small{7}};

\node at (10mm,26mm) {\small{2}};
\node at (10mm,10mm) {\small{6}};
\node at (10mm,2mm) {\small{5}};

\node at (14mm,26mm) {\small{2}};
\node at (14mm,10mm) {\small{6}};

\node at (18mm,10mm) {\small{\bf{3}}};

\node at (22mm,10mm) {\small{3}};
\node at (46mm,14mm){$\xlongrightarrow[j=6]{i=3,\,h=5}$};
\end{tikzpicture}&

\begin{tikzpicture}
\def\rectanglepath{-- +(4mm,0mm) -- +(4mm,4mm) -- +(0mm,4mm) -- cycle}
\draw [step=4mm,dotted] (0mm,0mm)grid (32mm,32mm);
\draw (0mm,24mm) \rectanglepath;
\draw (0mm,20mm) \rectanglepath;
\draw (0mm,12mm) \rectanglepath;
\draw (0mm,8mm) \rectanglepath;
\draw (0mm,4mm) \rectanglepath;
\draw (0mm,0mm) \rectanglepath;
\draw (4mm,24mm) \rectanglepath;
\draw (4mm,12mm) \rectanglepath;
\draw (4mm,8mm) \rectanglepath;
\draw (4mm,0mm) \rectanglepath;
\draw (8mm,24mm) \rectanglepath;
\draw (8mm,8mm) \rectanglepath;
\draw (8mm,0mm) \rectanglepath;
\draw (12mm,24mm) \rectanglepath;
\draw (12mm,8mm) \rectanglepath;
\draw (16mm,8mm) \rectanglepath;
\draw (20mm,8mm) \rectanglepath;

\node at (2mm,26mm) {\small{2}};
\node at (2mm,22mm) {\small{3}};
\node at (2mm,14mm) {\small{5}};
\node at (2mm,10mm) {\small{6}};
\node at (2mm,6mm) {\small{4}};
\node at (2mm,2mm) {\small{1}};

\node at (6mm,26mm) {\small{2}};
\node at (6mm,14mm) {\small{5}};
\node at (6mm,10mm) {\small{6}};
\node at (6mm,2mm) {\small{7}};

\node at (10mm,26mm) {\small{2}};
\node at (10mm,10mm) {\small{6}};
\node at (10mm,2mm) {\small{5}};

\node at (14mm,26mm) {\small{2}};
\node at (14mm,10mm) {\small{6}};

\node at (18mm,10mm) {\small{\bf{6}}};

\node at (22mm,10mm) {\small{3}};
\end{tikzpicture}
\end{tabular}
\caption{An example to illustrate Case 1.}\label{Case1}
\end{center}
\end{figure}
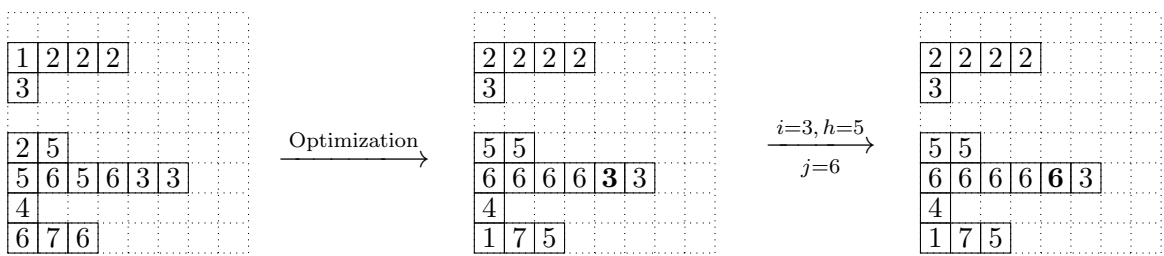

For example, Figure \ref{Case1} is an illustration of the construction of $F'$ in Case 1,
where the integer  in boldface signifies the integer $i$ that is changed to $j$.

Case 2.  $\beta_i>\beta_j$.  The filling $F'$ is obtained from $\opt(F)$ according to the following rule.
For $1\le m\le n$, let $F_m$ be the $m$-th column of $\opt(F)$, and $F'_m$  be the $m$-th column of $F'$.
\begin{itemize}
\item[(i)] If $F_m$ contains both  $i$ and $j$ or contains neither  $i$ nor $j$,
then let $F_m'=F_m$;

\item[(ii)] If $F_m$ contains only $i$, then $F_m'$ is obtained from $F_m$ by replacing $i$ with $j$;

\item[(iii)] If $F_m$ contains only $j$, then $F_m'$ is obtained from $F_m$ by  replacing $j$ with $i$.
\end{itemize}

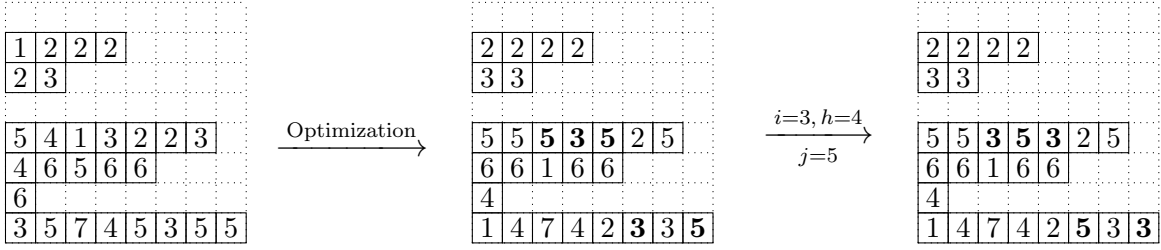
\begin{figure}[htbp]
\begin{center}
\begin{tabular}{ccc}
\begin{tikzpicture}
\def\rectanglepath{-- +(4mm,0mm) -- +(4mm,4mm) -- +(0mm,4mm) -- cycle}
\draw [step=4mm,dotted] (0mm,0mm)grid (32mm,32mm);
\draw (0mm,24mm) \rectanglepath;
\draw (0mm,20mm) \rectanglepath;
\draw (0mm,12mm) \rectanglepath;
\draw (0mm,8mm) \rectanglepath;
\draw (0mm,4mm) \rectanglepath;
\draw (0mm,0mm) \rectanglepath;
\draw (4mm,24mm) \rectanglepath;
\draw (4mm,20mm) \rectanglepath;
\draw (4mm,12mm) \rectanglepath;
\draw (4mm,8mm) \rectanglepath;
\draw (4mm,0mm) \rectanglepath;
\draw (8mm,24mm) \rectanglepath;
\draw (8mm,12mm) \rectanglepath;
\draw (8mm,8mm) \rectanglepath;
\draw (8mm,0mm) \rectanglepath;
\draw (12mm,24mm) \rectanglepath;
\draw (12mm,12mm) \rectanglepath;
\draw (12mm,8mm) \rectanglepath;
\draw (12mm,0mm) \rectanglepath;
\draw (16mm,12mm) \rectanglepath;
\draw (16mm,8mm) \rectanglepath;
\draw (16mm,0mm) \rectanglepath;
\draw (20mm,12mm) \rectanglepath;
\draw (20mm,0mm) \rectanglepath;
\draw (24mm,12mm) \rectanglepath;
\draw (24mm,0mm) \rectanglepath;
\draw (28mm,0mm) \rectanglepath;

\node at (2mm,26mm) {\small{1}};
\node at (2mm,22mm) {\small{2}};
\node at (2mm,14mm) {\small{5}};
\node at (2mm,10mm) {\small{4}};
\node at (2mm,6mm) {\small{6}};
\node at (2mm,2mm) {\small{3}};

\node at (6mm,26mm) {\small{2}};
\node at (6mm,22mm) {\small{3}};
\node at (6mm,14mm) {\small{4}};
\node at (6mm,10mm) {\small{6}};
\node at (6mm,2mm) {\small{5}};

\node at (10mm,26mm) {\small{2}};
\node at (10mm,14mm) {\small{1}};
\node at (10mm,10mm) {\small{5}};
\node at (10mm,2mm) {\small{7}};

\node at (14mm,26mm) {\small{2}};
\node at (14mm,14mm) {\small{3}};
\node at (14mm,10mm) {\small{6}};
\node at (14mm,2mm) {\small{4}};

\node at (18mm,14mm) {\small{2}};
\node at (18mm,10mm) {\small{6}};
\node at (18mm,2mm) {\small{5}};

\node at (22mm,14mm) {\small{2}};
\node at (22mm,2mm) {\small{3}};

\node at (26mm,14mm) {\small{3}};
\node at (26mm,2mm) {\small{5}};

\node at (30mm,2mm) {\small{5}};
\node at (46mm,14mm){$\xlongrightarrow[]{\text{Optimization}}$};
\end{tikzpicture}&
\begin{tikzpicture}
\def\rectanglepath{-- +(4mm,0mm) -- +(4mm,4mm) -- +(0mm,4mm) -- cycle}
\draw [step=4mm,dotted] (0mm,0mm)grid (32mm,32mm);
\draw (0mm,24mm) \rectanglepath;
\draw (0mm,20mm) \rectanglepath;
\draw (0mm,12mm) \rectanglepath;
\draw (0mm,8mm) \rectanglepath;
\draw (0mm,4mm) \rectanglepath;
\draw (0mm,0mm) \rectanglepath;
\draw (4mm,24mm) \rectanglepath;
\draw (4mm,20mm) \rectanglepath;
\draw (4mm,12mm) \rectanglepath;
\draw (4mm,8mm) \rectanglepath;
\draw (4mm,0mm) \rectanglepath;
\draw (8mm,24mm) \rectanglepath;
\draw (8mm,12mm) \rectanglepath;
\draw (8mm,8mm) \rectanglepath;
\draw (8mm,0mm) \rectanglepath;
\draw (12mm,24mm) \rectanglepath;
\draw (12mm,12mm) \rectanglepath;
\draw (12mm,8mm) \rectanglepath;
\draw (12mm,0mm) \rectanglepath;
\draw (16mm,12mm) \rectanglepath;
\draw (16mm,8mm) \rectanglepath;
\draw (16mm,0mm) \rectanglepath;
\draw (20mm,12mm) \rectanglepath;
\draw (20mm,0mm) \rectanglepath;
\draw (24mm,12mm) \rectanglepath;
\draw (24mm,0mm) \rectanglepath;
\draw (28mm,0mm) \rectanglepath;

\node at (2mm,26mm) {\small{2}};
\node at (2mm,22mm) {\small{3}};
\node at (2mm,14mm) {\small{5}};
\node at (2mm,10mm) {\small{6}};
\node at (2mm,6mm) {\small{4}};
\node at (2mm,2mm) {\small{1}};

\node at (6mm,26mm) {\small{2}};
\node at (6mm,22mm) {\small{3}};
\node at (6mm,14mm) {\small{5}};
\node at (6mm,10mm) {\small{6}};
\node at (6mm,2mm) {\small{4}};

\node at (10mm,26mm) {\small{2}};
\node at (10mm,14mm) {\small{\bf{5}}};
\node at (10mm,10mm) {\small{1}};
\node at (10mm,2mm) {\small{7}};

\node at (14mm,26mm) {\small{2}};
\node at (14mm,14mm) {\small{\bf{3}}};
\node at (14mm,10mm) {\small{6}};
\node at (14mm,2mm) {\small{4}};

\node at (18mm,14mm) {\small{\bf{5}}};
\node at (18mm,10mm) {\small{6}};
\node at (18mm,2mm) {\small{2}};

\node at (22mm,14mm) {\small{2}};
\node at (22mm,2mm) {\small{\bf{3}}};

\node at (26mm,14mm) {\small{5}};
\node at (26mm,2mm) {\small{3}};

\node at (30mm,2mm) {\small{\bf{5}}};

\node at (46mm,14mm){$\xlongrightarrow[j=5]{i=3,\,h=4}$};
\end{tikzpicture}&
\begin{tikzpicture}
\def\rectanglepath{-- +(4mm,0mm) -- +(4mm,4mm) -- +(0mm,4mm) -- cycle}
\draw [step=4mm,dotted] (0mm,0mm)grid (32mm,32mm);
\draw (0mm,24mm) \rectanglepath;
\draw (0mm,20mm) \rectanglepath;
\draw (0mm,12mm) \rectanglepath;
\draw (0mm,8mm) \rectanglepath;
\draw (0mm,4mm) \rectanglepath;
\draw (0mm,0mm) \rectanglepath;
\draw (4mm,24mm) \rectanglepath;
\draw (4mm,20mm) \rectanglepath;
\draw (4mm,12mm) \rectanglepath;
\draw (4mm,8mm) \rectanglepath;
\draw (4mm,0mm) \rectanglepath;
\draw (8mm,24mm) \rectanglepath;
\draw (8mm,12mm) \rectanglepath;
\draw (8mm,8mm) \rectanglepath;
\draw (8mm,0mm) \rectanglepath;
\draw (12mm,24mm) \rectanglepath;
\draw (12mm,12mm) \rectanglepath;
\draw (12mm,8mm) \rectanglepath;
\draw (12mm,0mm) \rectanglepath;
\draw (16mm,12mm) \rectanglepath;
\draw (16mm,8mm) \rectanglepath;
\draw (16mm,0mm) \rectanglepath;
\draw (20mm,12mm) \rectanglepath;
\draw (20mm,0mm) \rectanglepath;
\draw (24mm,12mm) \rectanglepath;
\draw (24mm,0mm) \rectanglepath;
\draw (28mm,0mm) \rectanglepath;

\node at (2mm,26mm) {\small{2}};
\node at (2mm,22mm) {\small{3}};
\node at (2mm,14mm) {\small{5}};
\node at (2mm,10mm) {\small{6}};
\node at (2mm,6mm) {\small{4}};
\node at (2mm,2mm) {\small{1}};

\node at (6mm,26mm) {\small{2}};
\node at (6mm,22mm) {\small{3}};
\node at (6mm,14mm) {\small{5}};
\node at (6mm,10mm) {\small{6}};
\node at (6mm,2mm) {\small{4}};

\node at (10mm,26mm) {\small{2}};
\node at (10mm,14mm) {\small{\bf{3}}};
\node at (10mm,10mm) {\small{1}};
\node at (10mm,2mm) {\small{7}};

\node at (14mm,26mm) {\small{2}};
\node at (14mm,14mm) {\small{\bf{5}}};
\node at (14mm,10mm) {\small{6}};
\node at (14mm,2mm) {\small{4}};

\node at (18mm,14mm) {\small{\bf{3}}};
\node at (18mm,10mm) {\small{6}};
\node at (18mm,2mm) {\small{2}};

\node at (22mm,14mm) {\small{2}};
\node at (22mm,2mm) {\small{\bf{5}}};

\node at (26mm,14mm) {\small{5}};
\node at (26mm,2mm) {\small{3}};

\node at (30mm,2mm) {\small{\bf{3}}};
\end{tikzpicture}
\end{tabular}
\caption{An example to illustrate Case 2.}\label{Case2}
\end{center}
\end{figure}

For example, Figure \ref{Case2} is an illustration of the construction of $F'$ in Case 2,
where the integers  in boldface signify the integers $i$ and $j$ that are interchanged.

It remains to show that $F'$ belongs to $\F(D(\alpha))$.  By the above constructions,
$F'_m$ is column strict. We need to show that
$F'_m$ satisfies the flag condition. This is evident if
$F'_m$ is obtained from $F_m$ by applying (i) or (iii).
We next prove the case when
$F'_m$ is obtained from $F_m$ by applying (ii).
Suppose that the entry $i$ in $F_m$ lies in the box $(l, m)$.

We first assert that
$\alpha_i<h$. Suppose otherwise that $\alpha_i\geq h$.
Then  the box $(i,h)$ belongs to $D(\alpha)$.
By the choice of the index $j$, the box $(i,h)$ is filled with $i$.
Since the box $(j,h)$ is also filled with $i$, the column $F_h$ of $F$
is not column strict, leading to a contradiction. This verifies the assertion.

Based on the above assertion, we claim that $l\geq j$.
 Suppose otherwise that $l<j$. Again, by the choice of the index $j$, we must have
 $l=i$. Since $(l,m)$ lies in the $i$-th row, we have $m\le\alpha_i$,
 which, together with  the above assertion that  $\alpha_i<h$, implies  $m< h$.
 Note that each box in row $j$ of $F$ that is to the left
 of $(j,h)$ is filled with $j$. This implies that the column $F_m$ contains both $i$ and $j$,
 which contradicts the fact that $F_m$ contains only $i$. Hence the claim is true.

By the above claim, we see that $F_m'$ satisfies the flag constraint, and so
$F'$ belongs to $\F(D(\alpha))$. Moreover, it is clear that $\beta=t_{i,j}(\wt(F'))$.
This completes the proof.
\qed

Based on Lemma \ref{LL}, we can now provide a proof of Theorem \ref{NN}.

\noindent
{\it Proof of Theorem \ref{NN}.}
Let $\beta$ be an exponent vector of $\kappa_\alpha(x)$. We aim to show that
$\beta\leq_\kappa \alpha$.
By Theorem \ref{CCC}, there exists a filling $F\in \F(D(\alpha))$ with
$\wt(F)=\beta$. The theorem is trivial in the case $\beta=\alpha$. We now
consider the case $\beta\neq \alpha$.

By Lemma \ref{LL},  there exists a
filling $F'$ in  $\F(D(\alpha))$   such that $\wt(F)=t_{i,j}(\wt(F'))$ or
$\wt(F)=m_{i,j}(\wt(F'))$. Notice that $\wt(F)>_{\mathrm{lex}} \wt(F')$,
where  $\geq_{\mathrm{lex}}$ is the lexicographic order on
 compositions. If $\wt(F')\neq \alpha$, then we can again invoke Lemma \ref{LL} to
 find a filling $F''\in \F(D(\alpha))$   such that $\wt(F')=t_{i,j}(\wt(F''))$ or
$\wt(F')=m_{i,j}(\wt(F''))$. Continuing this procedure, we can eventually
arrive at the (unique)  filling  $F_0\in \F(D(\alpha))$ with $\wt(F_0)=\alpha$.
Here, $F_0$ is the filling such that each box of $D(\alpha)$ is filled with its row index.
Hence $\beta=\wt(F)$ can be generated  from $\alpha=\wt(F_0)$ by applying a sequence of
moves $t_{i,j}$ and $m_{i,j}$, implying that $\beta\leq_\kappa \alpha$. This completes
the proof.
\qed



We  lastly  prove the sufficiency  of Theorem \ref{conj}.

\begin{theo}\label{SS}
If $\beta\leq_\kappa \alpha$, then $\beta$ is an exponent vector of
$\kappa_\alpha(x)$.

\end{theo}

\pf Let $\beta=\beta^{(0)}, \beta^{(1)},\ldots, \beta^{(k)}=\alpha$
be a sequence of compositions such that for $0\leq t\leq k-1$, $\beta^{(t)}$ is obtained from $\beta^{(t+1)}$
by applying $t_{i,j}$ or $m_{i,j}$.   The proof that $\beta$ is an exponent vector of
$\kappa_\alpha(x)$ is by induction on $k$. The case for $k=0$ (namely, $\beta=\alpha$)
is trivial.  We next consider the case $k>0$.

By hypothesis, $ \beta^{(1)}$ is an exponent vector of
$\kappa_\alpha(x)$.
It follows from  Theorem \ref{CCC} that there is a filling $F$ in  $\F(D(\alpha))$
such that $\wt(F)=\beta^{(1)}$. To complete proof, it suffices to
construct a filling $F'$ in  $\F(D(\alpha))$
such that $\wt(F')=\beta$. Write $\beta^{(1)}=(v_1,\ldots, v_n)$.  We have two cases.

Case 1. $\beta=t_{i,j}(\beta^{(1)})$. In this case, $v_i<v_j$. Suppose that there are
$c_1$ columns of $F$ that contain   $i$ but do not contain $j$. Then
there are
$c_1+v_j-v_i$ columns of $F$ that contain   $j$ but do not contain $i$.
Choose any $v_j-v_i$ such columns, and let $F'$ be obtained from
$F$ by replacing   $j$ in each of these columns with $i$. It is readily
seen that $F$ belongs to $\F(D(\alpha))$ with $\wt(F')=\beta$.

Case 2. $\beta=m_{i,j}(\beta^{(1)})$. In this case, $v_i<v_j-1$. Hence there is at least one
column of $F$ that contains   $j$ but does not contain $i$. Choose one such column, and
let $F'$ be obtained from
$F$ by replacing   $j$ in this  column  with $i$. It is also clear that
 $F$ belongs to $\F(D(\alpha))$ with $\wt(F')=\beta$.
This completes the proof.
\qed

\vspace{.2cm} \noindent{\bf Acknowledgments.}
This work was
supported by the 973 Project  and the National Science Foundation of China (Grant No. 11971250).

\footnotesize{
N{\scriptsize EIL} J.Y. F{\scriptsize AN}, D{\scriptsize EPARTMENT OF} M{\scriptsize ATHEMATICS}, S{\scriptsize ICHUAN} U{\scriptsize NIVERSITY}, C{\scriptsize HENGDU} 610064, P.R. C{\scriptsize HINA.} Email address: fan@scu.edu.cn.

P{\scriptsize ETER} L. G{\scriptsize UO}, C{\scriptsize ENTER FOR} C{\scriptsize OMBINATORICS}, N{\scriptsize ANKAI} U{\scriptsize NIVERSITY}, T{\scriptsize IANJIN} 300071, P.R. C{\scriptsize HINA.}  Email address: lguo@nankai.edu.cn.
}

S{\scriptsize IMON} C.Y. P{\scriptsize ENG}, C{\scriptsize ENTER FOR} A{\scriptsize PPLIED} M{\scriptsize ATHEMATICS}, T{\scriptsize IANJIN} U{\scriptsize NIVERSITY}, T{\scriptsize IANJIN} 300072, P.R. C{\scriptsize HINA.}  Email address: pcy@tju.edu.cn.

S{\scriptsize OPHIE} C.C. S{\scriptsize UN}, D{\scriptsize EPARTMENT OF} M{\scriptsize ATHEMATICS}, U{\scriptsize NIVERSITY OF} F{\scriptsize INANCE AND} E{\scriptsize CONOMICS}, T{\scriptsize IANJIN} 300222, P.R. C{\scriptsize HINA.}  Email address: suncongcong@mail.nankai.edu.cn.

\end{document}